\documentclass[conference]{IEEEtran}
\ifCLASSINFOpdf
\else
\fi
\usepackage{hyperref}
\usepackage{amsmath}
\usepackage{amssymb}
\usepackage{bm}
\usepackage{comment}
\usepackage{graphicx}
\usepackage{algorithm}
\usepackage{algorithmic}
\newtheorem{theorem}{Theorem}[section]

\newtheorem{lemma}[theorem]{Lemma}

\DeclareMathOperator*{\argmin}{arg\,min}
\begin{document}
%
\title{Stochastic Decision-Making Model for Aggregation of Residential Units with PV-Systems and Storages}

\author{\IEEEauthorblockN{Hossein Khazaei}
\IEEEauthorblockA{\textit{Dep. of Electrical and Computer Eng.} \\
	\textit{Stony Brook University}\\
	Stony Brook, NY, USA \\
	\href{mailto:me@somewhere.com}{hossein.khazaei@stonybrook.edu}}
\and
\IEEEauthorblockN{Ramin Moslemi}
\IEEEauthorblockA{\textit{Energy Management Department} \\
	\textit{NEC Laboratories America, Inc.}\\
	San Jose, CA, USA \\
	\href{mailto:rmoslemi@nec-labs.com}{rmoslemi@nec-labs.com}}
\and
\IEEEauthorblockN{Ratnesh Sharma}
\IEEEauthorblockA{\textit{Energy Management Department} \\
	\textit{NEC Laboratories America, Inc.}\\
	San Jose, CA, USA \\
	\href{mailto:ratnesh@nec-labs.com}{ratnesh@nec-labs.com}}
}

\maketitle
\thispagestyle{plain}
\pagestyle{plain}

\begin{abstract}
Many residential energy consumers have installed photovoltaic (PV) panels and energy storage systems. These residential users can aggregate and participate in the energy markets. A stochastic decision making model for an aggregation of these residential units for participation in two-settlement markets is proposed in this paper. Scenarios are generated using Seasonal Autoregressive Integrated Moving Average (SARIMA) model and joint probability distribution function of the forecast errors to model the uncertainties of the real-time prices, PV generations and demands.
The proposed scenario generation model of this paper treats forecast errors as random variable, which allows to reflect new information observed in the real-time market into scenario generation process without retraining SARIMA or re-fitting probability distribution functions over the forecast errors. This approach significantly improves the computational time of the proposed model.
A simulation study is conducted for an aggregation of 6 residential units, and the results highlights the benefits of aggregation as well as the proposed stochastic decision-making model.
\end{abstract}


%
\IEEEpeerreviewmaketitle

\section{Introduction} \label{sec::Intro}
Technological developement in the past decade made rooftop PV systems and battery energy storage systems (BESS) economically attractive to the residential units. These PV systems are usually coupled with a BESS system to increase their economical efficiency and flexibility. They can aggregate and participate in the energy markets, and the aggregation is usually represented by an aggregator in the markets.

When the number of residential households with BESS and PV system is large enough, they can aggregate and directly participate in the energy markets, instead of buying energy from utilities, which decreases their costs. Another benefit of aggregation is that the aggregation of random variables decreases their total volatility \cite{Bitar2012, Khazaei2019}. The sources of volatilities for the residential households are: uncertainty in the prices, PV generations, and the demand for the coming hours. These parameters depend on exogenic forces, especially weather. Now the question is: how an aggregator can efficiently represent the households in the markets considering these uncertainties?

In \cite{PVandBESS_Paper2} a decision making model for the aggregation of residential units with BESS and PV systems is proposed, where the forecast errors of the PV generation and the electricity prices are handled using a model predictive model (MPC). MPC demands high computational power and solving the MPC is challenging when there are uncertainties. A different approach to handle the uncertainties is to use robust optimization, which uses the notion of \textit{uncertainty set} to aggregate the adverse events that we are particularly interested to consider.
A robust optimization framework for real-time operation of a residential battery storage unit in the presence of price uncertainties is proposed in \cite{PVandBESS_Paper4}.
Another robust decision making model for an aggregation of BESS and PV systems for participating in the DA market is proposed in \cite{PVandBESS_Paper3}, where the uncertainties in the DA prices, PV generations and loads are modeled using scenarios.
A robust optimization approach is proposed in \cite{PVandBESS_Paper6} to derive the optimal behavior of a distributed energy resource (e.g. PV systems) in the DA market where the market prices and PV generations are uncertain.
The solution obtained from the robust optimization usually is less optimal than the solution derived from the stochastic optimization because stochastic optimization benefits from using the characterization of the distribution of the uncertainty \cite{Ackooij18}.
\\
Another approach for managing the uncertainties is to use \textit{chance-constraint optimization}, which uses the feasibility probability of the solution as a parameter to setting the trade-off between the robustness of the solution and its optimality \cite{Bruninx18}. Similar to the robust optimization, the solution of the chance-constraint optimization is often less optimal than that of the stochastic optimization.
\subsection{Contributions}
This paper presents a stochastic decision-making model for an aggregation of BESS and PV systems to participate in the two-settlement market considering the cost of using the network facilities and the battery degradation cost. The stochastic optimization model is less conservative than the worst-case-oriented robust optimization models or the chance-constraint optimization models \cite{RO_BenTal}.
The Seasonal Autoregressive Integrated Moving Average (SARIMA) model is trained using historical data for forecasting RT market prices, PV generations and loads. Forecast errors are represented by random variables, and the joint probability distribution function of the forecast error is derived to generate a set of volatilities. Such volatilities is then added to the forecasts to create scenarios for the stochastic decision making model. By modeling the forecast errors as random variables the proposed model incorporates new information of the RT market about real-time prices, PV generations and loads into scenario generation process without retraining the SARIMA model. Such approach generates scenarios that are consistent with observations of prior hours in the RT market in a computationally efficient manner. A modified Simultaneous Backward Scenario Reduction algorithm is used to decrease the number of scenarios while keeping the set of preserved scenarios relatively accurate.

The proposed model of this paper takes advantages of stochastic optimization approaches while efficiently addressing the key limitation of these models, which is the computational complexity in handling large number of scenarios.
\section{System Model} \label{sec::SysMdl}
\subsection{Two-Settlement Market Model}  \label{subsec::TwoSetMark}
We consider a two-settlement market model consisting of a Day-Ahead (DA) market and a Real-Time (RT) market. DA market is a forward market where generation units and loads bid into it and the Independent System Operator (ISO) clears the market. The results of the DA market are promises that must be delivered next day at time of delivery. Any imbalance between the DA promise and the actual performance at the time of delivery is cleared in the RT market.

As a benchmark, we first consider a single residential unit $i$ participating in this two-settlement market. This residential unit owns a BESS and PV system. It also has a random load that needs to met. The two-settlement market is detailed in the following. \\
\textbf{At DA market:} The residential unit $i$ observes the DA prices $\left\{p^D_1,\cdots, p^D_{24}\right\}$, and bids for selling/buying power at this market to be delivered/received at the corresponding hours of the next day. The bid for hour $t$ is shown with $c_{i,t}$. If $c_{i,t} > 0$ ($c_{i,t} < 0$), it means the residential unit bids for selling (buying) $|c_{i,t}|$ unit of power at the DA market to be delivered (received) at hour $t$ of the next day (\textit{i.e.} at delivery time). 
\\
\textbf{At RT market, hour} $t$\textbf{:} For this hour, the residential unit observes the realizations of its demand $\left(d_{i,t}\right)$, PV generation $\left(v_{i,t}\right)$, and the RT price $\left(p^R_t\right)$. At this stage, the residential unit has a chance to participate in the RT market and buy/sell power. At this market, the household submits a bid, $x_{i,t}$ to the ISO for which positive values of $x_{i,t}$ corresponds to selling $x_{i,t}$ unit of power at RT market and negative values of $x_{i,t}$ correspond to buying $|x_{i,t}|$ units of power from the RT market. \\
The key equation stating the relationship between DA commitment, RT bid, PV generation, and BESS for hour $t$ is in the following:
\begin{align} \label{eq::PwrBlnc}
s_{i,t} + \eta_i \cdot \left(- c_{i,t} + v_{i,t} - d_{i,t} - x_{i,t}\right) = s_{i,t+1}
\end{align}
where $s_{i,t}$ is the energy stored in the BESS at the begining of hour $t$. Note that $s_{i,t} - s_{i,t+1}$ represents the battery discharge at hour $t$ and there is an underlying degradation cost associated with it. There is also an upper and lower level limit $\left(i.e. \ s_i^{min}, s_i^{max}\right)$ to the energy stored at the BESS. $\eta_i$ is the charging/discharging efficiency rate for the BESS. 
\\
\textbf{Assumptions:} In this paper we assume that the residential unit is price taker, meaning that it does not have market power to affect the market prices, and the ISO allocates the residential unit whatever it bids in both DA and RT markets. Later when considering the aggregation of the residential units, we do not consider the transmission constraints between the locations of the residential units. 

Now we consider an aggregator representing a set of  residential units in the two-settlement market.

\subsection{Aggregator Participating In Two-Settlement Market}  \label{subsec::AggParticMarket}
It is shown that an aggregation of individual units with underlying volatilities is relatively less volatile comparing to the sum of those individual units \cite{Bitar2011}. 
In this paper we implement this idea to improve the performance of a set of residential units with BESS and PV systems participating in the DA-RT market. 
The setting of the aggregation is as follows: \\
A set of $N$  residential units with BESS and PV systems are represented by an aggregator in the DA-RT market. The BESS and the PV system of the aggregator is the aggregation of the corresponding systems of the residential units,
\begin{align}
s_{A}^{min} = \sum_{i = 1}^{N} s_i^{min}, \hspace{10pt} &  s_{A}^{max} = \sum_{i = 1}^{N} s_i^{max}, \nonumber \\
v_{A,t} = \sum_{i = 1}^{N} v_{i,t}, \hspace{10pt} & d_{A,t} = \sum_{i = 1}^{N} d_{i,t}
\end{align}
In the next section we detail the stochastic decision making model used by the aggregator in the DA-RT market.

\section{Stochastic Decision Making Model} \label{sec::StocDecMakModel}
For a 24-hour day, the aggregator solves 25 stochastic decision-making problems at 25 different time slots, one at the DA market and 24 for each hour at the RT market. Decision variables of the aggregator and the uncertain parameters it faces at each decision-making problem are listed in Table \ref{tab:tab1}.
\begin{table}[htbp]
	\caption{Parameters of Different Decision-Making Problems}
	\vspace{-9pt}
	\begin{center}
		\centering
		\begin{tabular}{|c|c|c|c|}
			\hline
			& Known Key & Unknown Key & Decision \\
			& Parameters & Parameters & Variables \\
			\hline
			DA & & $v_{A,1}, \cdots, v_{A,24}$ & \\
			Market & $p^D_1, \cdots, p^D_{24}$ & $d_{A,1},\cdots, d_{A,24}$ &  $c_{A,1}, \cdots, c_{A,24}$ \\
			\hline
			RT & $v_{A,t}, d_{A,t}$ & $v_{A,t+1}, \cdots, v_{A,24}$ & \\
			Market & $c_{A,t}, \cdots, c_{A,24}$ & $d_{A,t+1},\cdots, d_{A,24}$ &  $x_{A,t}$ \\
			$\left(\mbox{hour } t\right)$ & $s_{A,t}, p^R_t$ & $p^R_{t+1}, \cdots, p^R_{24}$ & \\
			\hline
		\end{tabular}
		\label{tab:tab1}
	\end{center}
\end{table}
At the DA market the aggregator tries to minimize its expected cost for the next day given the constraints that it faces, especially regarding the energy storage. 
In this paper we use scenario generation to numerically solve the stochastic decision-making model. The scenario generation model generates different scenarios with different probabilities to reflect the randomness of the uncertain parameters.
\\
At the DA market, agregator solves the following problem:
\begingroup
\allowdisplaybreaks
\begin{align} \label{eq::DAdecMakMdl}
&\min_{c_1,\cdots,c_{24}} \pi_A = -\sum_{i \in T} p^D_t  c_{A,i}  + \hspace{-5pt} \sum_{k \in \Omega^D} \hspace{-4pt} \epsilon^k \cdot \Big( \hspace{-3pt} - \sum_{i \in T} p^{R,k}_i \cdot x_{A,i}^k   \nonumber  \\
& \hspace{40pt} + \alpha \cdot \sum_{i \in T} \left( s_{A,i+1}^k - s_{A,i}^k \right)^2 + \beta \cdot \sum_{i \in T} \left| c_{A,i} + x_{A,i} \right| \hspace{-3pt} \Big) \nonumber \\
& \mbox{s.t.} \nonumber \\
& s_{A,i}^k + \eta \cdot \left(- c_{A,i} + v_{A,i}^k - d_{A,i}^k - x_{A,i}^k\right) = s_{A,i+1}^k  \nonumber \\
& \hspace{163pt}  \ \ \ \forall k \in \Omega^D , i \in T_0 \nonumber \\
& s_{A}^{min} \leq s_{A,i}^k \leq s_A^{max} \hspace{90pt} \forall k \in \Omega^D , i \in T_0 \nonumber \\
& | c_{A,i} | \leq c_{A}^{max} \hspace{119pt} \forall i \in T_0 
\end{align}
\endgroup
$\pi_A$ is the estimated expected cost of the aggregator. We use the term \textit{estimated} to emphasize that here we are estimating the expected cost using a set of scenarios ($\Omega^D$) modeling the uncertain parameters. $s^k_{25}$ is the net energy at the battery at the end of the day under scenario $k$.
$\epsilon^k$ is the probability of scenario $k$. $\alpha$ and $\beta$ are the cost multipliers corresponding to the battery degradation and distribution network usage.

At hour $t$ of the RT market, aggregator observes the PV generation, demand and the RT price for that hour along with the net energy of BESS at the beginning of hour $t$. The decision-making problem that the aggregator is facing is as follows:
\begingroup
\allowdisplaybreaks
\begin{align} \label{eq::RTdecMakMdl}
&\min_{x_{A,t}} \ \ - p^R_t  x_{A,t}  + \sum_{k \in \Omega^R_k} \epsilon^k \cdot \Big( - \sum_{i \in T_t} p^{R,k}_i  x_{A,i}^k   \nonumber  \\
& \hspace{28pt} + \alpha \cdot \sum_{i \in T_t} \left( s_{A,i+1}^k - s_{A,i}^k \right)^2 + \beta \cdot \sum_{i \in T_t} \left|c_{A,i} + x_{A,i}\right| \Big) \nonumber \\
& \mbox{s.t.} \nonumber \\
& s_{A,t} - c_{A,t} + v_{A,t} - d_{A,t} - x_{A,t} = s_{A,t+1}^k, \ \ \ \forall k \in \Omega^D  \nonumber \\
& s_{A,i}^k - c_{A,i} + v_{A,i}^k - d_{A,i}^k - x_{A,i}^k = s_{A,i+1}^k, \hspace{-2pt} \ \forall k \in \Omega^D , i \in T_t \nonumber \\
& s_{A}^{min} \leq s_{A,i}^k \leq s_A^{max} \hspace{60pt} \forall k \in \Omega^D , i \in T_t 
\end{align}
\endgroup
Where $T_t \triangleq \left\{t+1,\cdots,24\right\}$ is the set of remaining hours of the day. In order to use this model we need to generate a set of scenarios for the uncertain parameters.
%
 %
\subsection{Scenario Generation} \label{subsec::ScenGener}
In both the DA and RT decision-making problems there are three sources of uncertainties: RT prices, PV generations, and the loads, all for the coming hours of the RT market. Note that we do not consider the uncertainty in the DA prices, but the results of this paper can be extended to include uncertainty in the DA prices.
In this paper we assume that these random variables are independent. 
Steps for the scenario generation for each random variable is detailed in Algorithm \ref{AlgScenDA}.
First we train the Seasonal Autoregressive Integrated Moving Average (SARIMA) model \cite{PVandBESS_Paper5} and derive the joint probability distribution function (\textit{pdf}) of the forecast error of the SARIMA model for the 24 hours. 
The trained SARIMA generates a forecast for the time series for the next day, and noises based on the joint \textit{pdf} of the forecast error is added to it to generate the scenarios \cite{ConejoBook}.
%
\begin{algorithm}
\begin{algorithmic}[1]
\caption{Scenario generation in DA market} 
\label{AlgScenDA}
 \renewcommand{\algorithmicrequire}{\textbf{Input:}}
 \renewcommand{\algorithmicensure}{\textbf{Output:}}
 \REQUIRE Historical Data of PV generations, Demand, and RT Prices.
 \ENSURE  Set of scenarios
 \\ \textit{DA Forecast}: 
 \STATE Use the historical data of the uncertain parameter to train the SARIMA.
 \STATE Calculate the forecast error of the historical data and fit a multivariate normal distribution function to the forecast errors for the 24 hours.
 \\ \textit{Scenario Generation}:
 \STATE Use SARIMA make a forecast for the uncertain parameter for the coming hours.
 \STATE Generate a set of scenarios by sampling from the \textit{pdf} of the forecast error and add it to the forecast derived from SARIMA model.
 \RETURN Set of scenarios for the next day. 
 \end{algorithmic} 
 \end{algorithm}
 %
 %
 
 The main contribution of this paper is the model we propose for scenario generation at each hour of RT market. First, the forecast errors of the SARIMA for each hour is represented by a random variable. At hour $t$ of the RT market the aggregator observes the \textit{realized} values of the forecast errors for current and previous hours, \textit{i.e.} hours $1$ to $t$. Using such observations the aggregator derives the updated \textit{pdf} of the forecast errors for hours $t+1$ to $24$ \textit{conditioned} on prior observations. Second, it is assumed that forecast errors for the 24 hours has a \textit{normal pdf}. Such assumption allows to derive the \textit{conditional pdf} of forecast errors using closed form formula provided in Lemma \ref{lemma::PdfUpdate} without re-fitting the \textit{pdf}. Such modeling approach allows us to make a forecast in the DA market using SARIMA, and reflect all new information of the RT market into the scenario generation model without retraining SARIMA or re-fitting the probability distribution functions.
 %
 %
 Steps for scenario generation in RT market at hour $t$ is detailed in Algorithm \ref{AlgScenRT}.
\begin{algorithm}
\begin{algorithmic}[1]
\caption{Scenario generation at time $t$ of RT market} 
\label{AlgScenRT}
 \renewcommand{\algorithmicrequire}{\textbf{Input:}}
 \renewcommand{\algorithmicensure}{\textbf{Output:}}
 \REQUIRE \begin{itemize}
     \item[$ $]
     \item DA forecasts of SARIMA for the PV generation, demand, and RT prices.
     \item Realizations of PV generation, demand, and RT prices for hours $1,\cdots,t$.
 \end{itemize}
 \ENSURE  Set of scenarios
 \STATE Using the realizations of PV generation, demand, and RT prices for hours $1,\cdots,t$ calculate the realized values of forecast error for such hours.
 \STATE Calculate the conditional joint \textit{pdf} of the forecast errors of the PV generation, demand, and RT prices for hours $t+1,\cdots,24$ using Lemma \ref{lemma::PdfUpdate}.
 \STATE Generate set of scenarios by sampling from the pdf of the forecast error and add it to the forecast derived from SARIMA model.
 \RETURN Set of optimal scenarios for random variables for hours $t+1, \cdots, 24$ of RT market. 
 \end{algorithmic} 
 \end{algorithm}
\begin{lemma} \label{lemma::PdfUpdate}
	If $X = \begin{bmatrix}
	X_{1} \\
	X_{2} \\
	\end{bmatrix}$, and $X$ has a normal distribution $X \sim \mathcal{N} \left(\mu_X, \Sigma_X\right)$ where $\mu_X = \begin{bmatrix}
	\mu_{1} \\
	\mu_{2} \\
	\end{bmatrix}$ and $\Sigma_X = \begin{bmatrix} 
	\Sigma_{11} & \Sigma_{12} \\
	\Sigma_{21} & \Sigma_{22} 
	\end{bmatrix}$
	then the conditional pdf of $X_2$ given the realization of $X_1$ is :
	\begin{align}
	\mbox{P} \left( X_2 \big| X_1 = \beta \right) = \mathcal{N} \big(\hat{\mu}_{2} , \hat{\Sigma}_{2}\big),
	\end{align}
	where
	\begingroup
    \allowdisplaybreaks
	\begin{gather}
	\hat{\mu}_{2} = \mu_2 + \Sigma_{21} \times \Sigma_{11}^{-1} \times \left(\beta - \mu_1\right) \nonumber \\
	\hat{\Sigma}_{2} = \Sigma_{22} + \Sigma_{21} \times \Sigma_{11}^{-1} \times \Sigma_{12}
	\end{gather}
	\endgroup
\end{lemma}
%
%
%
\subsection{Scenario Reduction} \label{Sec::ScenRed}
The algorithm we use for the scenario reduction is a modified version of \textit{Simultaneous Backward Reduction} detailed in \cite{ScenarioReductionRef}. We first define the \textit{distance} between two scenarios generated from a stochastic process. Let $\zeta_{t\rightarrow{}T}^i$ and $\Tilde{\zeta}_{t\rightarrow{}T}^j$ be two scenarios of $n$ dimensional random variable from time $t$ to $T$. The \textit{distance} of these two scenarios is defined as:
\begin{align}
    c_{t\rightarrow{}T} \big( \zeta_{t\rightarrow{}T}^i , \Tilde{\zeta}_{t\rightarrow{}T}^j \big) \triangleq \sum_{\tau = t}^{T} \big( \zeta_{\tau}^i - \Tilde{\zeta}_{\tau}^j  \big)^2
\end{align}
We use such distance for scenario reduction. The steps for scenario reduction is detailed in Algorithm \ref{AlgScenRed}.
\begin{algorithm}
\begin{algorithmic}[1]
\caption{Simultaneous Backward Scenario Reduction} 
\label{AlgScenRed}
 \renewcommand{\algorithmicrequire}{\textbf{Input:}}
 \renewcommand{\algorithmicensure}{\textbf{Output:}}
 \REQUIRE Set of scenarios, \textit{i.e.} $S_S$, including $K_S$ scenarios for a stochastic process from hour $t$ to hour $24$, and their probabilities, \textit{i.e.} $\{p_i\}$.
 \ENSURE  Set of $K_P$ preserved scenarios. \\
 \textit{Initialization}: $S_D = \emptyset$ (set of deleted scenarios)
 \STATE Calculate the distance of pair scenarios in $S_S$.
 \[
 c_{ij} \triangleq c_{t\rightarrow{}24} \big( \zeta_{t\rightarrow{}24}^i , \Tilde{\zeta}_{t\rightarrow{}24}^j \big), i,j \in S_S
 \]
 \vspace{-8pt}
 \FOR {$k = 1$ to $K_P$:}
 \STATE Compute the distance between the pairs in the set of deleted scenarios and the set of remaining scenarios: \label{AlgScenRed_Dist}
 \[
 c_{ij}^{[k]} \triangleq \min_{u \notin S_D \cup \{j\}} c_{iu}, \forall i \in S_D \cup \{j\}, \forall j \notin S_D.
 \]
 \STATE Compute $z_j^{[k]} \triangleq \sum_{i \in S_D \cup \{j\}} p_i c_{ij}^{[k]}, j \notin S_D$
 \vspace{4pt}
 \STATE Select $j_i \in \argmin_{j \notin S_D}{z_j^{[k]}}$.
 \STATE Update the set of deleted scenarios: $S_D \xleftarrow{} S_D \cup \{j_i\}$
 \ENDFOR
 \STATE The set of reduced scenarios is: $S_P = S_S - S_D$.
 \STATE Calculate the probabilities for the set of preserved scenarios: \label{AlgScenRed_UpdProb}
 \[
 \hat{p}_i = p_i + \sum_{u \in J_{S,i}}{p_i}, \ \ \forall i \in S_P
 \]
 \vspace{-3pt}
 where
 \[
 J_{S,i} \triangleq \big\{ v \in S_D \ \big| \ i \in \argmin_{u \notin S_D}{c_{t\rightarrow{}T} \big( \zeta_{t\rightarrow{}T}^u , \Tilde{\zeta}_{t\rightarrow{}T}^v \big)}      \big\}
 \]
 \vspace{-5pt}
 \RETURN $S_P$ and $\{\hat{p}_i\}$.
 \end{algorithmic} 
 \end{algorithm}
%
At each iteration the algorithm deletes the scenario that has the lowest \textit{value} to the remaining scenarios. \textit{Value} of a scenario, defined in step \ref{AlgScenRed_Dist} of Algorithm \ref{AlgScenRed}, is a combination of its distance from deleted scenarios and their probabilities. At the end of each iteration, the scenario with lowest \textit{value} is added to the set of deleted scenarios. The probability of the preserved scenarios is updated in step \ref{AlgScenRed_UpdProb} where the probability of each deleted scenario is added to scenario in the set of preserved scenarios with lowest distance from it.
\section{Simulations}  \label{sec::Simul}
We perform the simulation using the data of PV generations of 6 households at California electricity market for the year 2012 \cite{PricesData}. Each residential unit owns 1KW-peak rooftop PV system and a BESS. We assume that the characteristics of the battery storage systems are similar for all the 6 households. Table \ref{tab:tab2} summarizes the properties of the battery storage systems. 
\begin{table}[htbp]
\vspace{-4pt}
	\caption{Parameters of Each Battery Storage Systems}
	\vspace{-9pt}
	\begin{center}
		\centering
		\begin{tabular}{|l|c|}
			\hline
			\hspace{70pt} Parameter & Value  \\
			\hline
			Maximum Battery capacity (kWh) & 5 \\
			\hline
			Charging efficiency (\%) & 90 \\
			\hline
			Minimum and maximum storage & 10-90  \\
			\hline
			Quadratic cost coefficient ($\alpha$) (cf. \eqref{eq::DAdecMakMdl}, \eqref{eq::RTdecMakMdl}) & 0.4 \\
			\hline
		\end{tabular}
		\label{tab:tab2}
	\end{center}
	\vspace{-8pt}
\end{table}
We also assume that the cost coefficient for using distribution network is $\beta = 0.05$ (cf. \eqref{eq::DAdecMakMdl}, \eqref{eq::RTdecMakMdl}).
We trained the SARIMA model using historical data of PV generations, RT prices and residential consumptions for $6$ months. The stochastic decision making problem is solved for Feb. 2012. 
At each step we generate 50 scenarios for each uncertain variable, and reduce it to 5, hence the total number of scenarios at each step is 125.
\begin{figure}[htbp]
	\centerline{\includegraphics[scale=0.53]{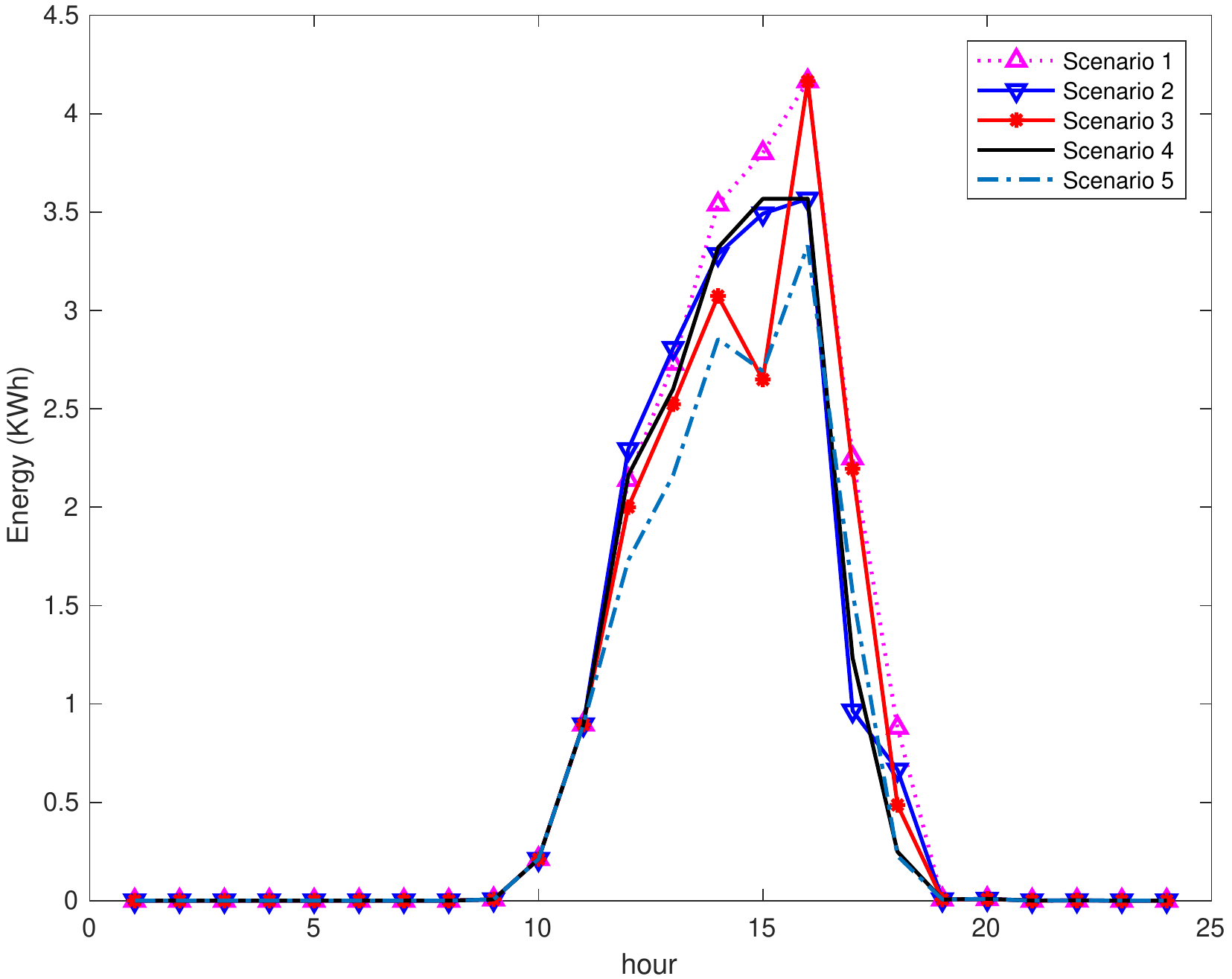}}
	\caption{Scenarios for the PV generation.}
	\label{fig3}
\end{figure}
Figure \ref{fig3} illustrates 5 scenarios for the PV generation for a particular day. These scenarios are derived after generating 50 scenarios and reducing them to 5 using Simultaneous Backward Scenario Reduction method. 
Note that the probability of each scenario may be different. 
\\
In Figure \ref{fig1} the performance of the aggregator for one day is shown. 
At the beginning of the day, when the prices are low, aggregator tries to charge the BESS. At the middle of the day aggregator uses the power produced by the PV systems, and tries not use the BESS. In the afternoon aggregator starts using the energy stored in the BESS as the market prices are high and the power generated by the PV systems has already started falling.
\begin{figure}[htbp]
	\centerline{\includegraphics[scale=0.45]{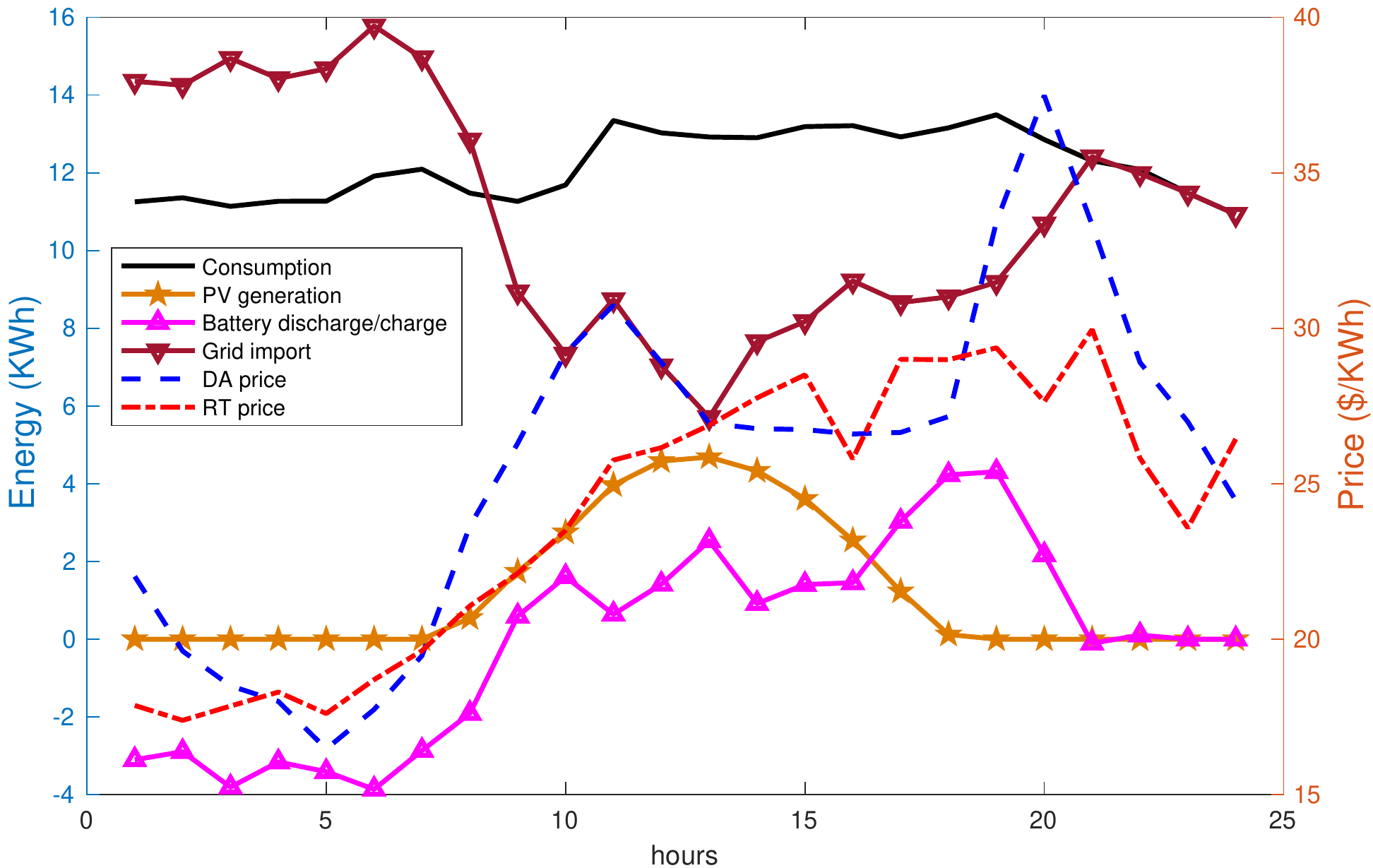}}
	\caption{Performance of aggregator for one day.}
	\label{fig1}
\end{figure}
%
%
%
Table \ref{tab:tab3} compares the total cost of supplying power to the individual units in three cases:
\\
\textbf{Case 1 (the proposed setting):} Aggregator supplies power to the individual units and uses the proposed stochastic decision-making model.
\\
\textbf{Case 2:} There is an aggregator which supplies power to the individual units. However, there is no stochastic decision-making model. Instead, there is a deterministic decision-making model where it is \textit{assumed} that the realization for the uncertain parameters in the past 24 hours will be repeated for the coming hours.
\\
\textbf{Case 3:} There is no aggregation, and each individual unit uses the proposed stochastic decision-making model.
\\
The total cost for these cases is shown in Table \ref{tab:tab3}.
%
%
\begin{table}[htbp]
\vspace{-4pt}
	\caption{Total cost of aggregator for different settings}
	\vspace{-9pt}
	\begin{center}
		\centering
		\begin{tabular}{|l|c|}
			\hline
			  & Cost ($\$$)  \\
			\hline
			Case 1: Aggregation $+$ Stochastic Decision Making Model & 169100 \\
			\hline
			Case 2: Aggregation, No Stochastic Decision Making Model & 227300 \\
			\hline
			Case 3: No Aggregation, Stochastic Decision Making Model & 185400  \\
			\hline
		\end{tabular}
		\label{tab:tab3}
	\end{center}
	\vspace{-8pt}
\end{table}
Note that the total cost for case 1 is lower than case 2 and 3, and the reason is that in case 1, we are benefiting from both the aggregation and using the proposed stochastic decision-making model. By comparing case 1 and case 2 we can see the benefit of using the proposed stochastic decision-making model. Also, by comparing case 1 and case 3 we can see the benefit of  aggregation. Note that by comparing cases 2 and 3 we can see that the total cost at case 3 is lower than the total cost for case 2. This implies that the benefit that the individual units receive from proper modeling the uncertain parameters by using the proposed stochastic decision-making model is higher than the benefit they obtain by aggregation. 			
%
%
%
%

\section{Conclusions}   \label{sec::Conclu}
A stochastic decision making model for an aggregation of residential units with Battery Energy Storage Systems (BESS) and Photovoltaic (PV) systems to participate in the two-settlement markets is presented. The uncertainties of the PV generations, real-time (RT) prices, and the demand is modeled using scenarios. Seasonal Autoregressive Integrated Moving Average Model (SARIMA) is used to derive a time-series forecast of the unknown parameters for the coming hours in the RT market, and then multiple scenarios are generated using the multivariate probability distribution function (pdf) of the forecast errors. With new information in the RT market, this joint \textit{pdf} is updated. This ensures the new information is incorporated in the model, while there is no need to retraining of SARIMA or re-fitting the pdf to the forecast errors. Also a modified Simultaneous Backward Scenario Reduction algorithm is proposed to decrease the number of scenarios so that the stochastic optimization problem becomes computationally feasible.
%
%
%
%
%
%
%
%
%
%
%
\bibliographystyle{IEEEtran}
{\bibliography{TPS}}

\end{document}